\documentclass{amsart}

\usepackage{amsmath, amsthm}
\usepackage{amssymb, latexsym}
\usepackage{amsfonts}

\newtheorem{theorem}[equation]{Theorem}

\newtheorem{corollary}[equation]{Corollary}

\theoremstyle{definition}
\newtheorem{definition}[equation]{Definition}

\newcommand{\thmref}[1]{\rm Theorem~\ref{#1}}
\newcommand{\defref}[1]{Def\-i\-ni\-tion \ref{#1}}

\newcommand{\corref}[1]{Corollary~\ref{#1}}

\newcommand{\beql}[1]{\begin{equation}\label{#1}}
\newcommand{\eeq} {\end{equation}}

    \font\Aaa=msam10

\font\Bbb=msbm10

\newcommand\Z{\hbox{\Bbb Z}}

\font\Aaa=msam10

\def\qed{\hbox{~~\Aaa\char'003}}

\font\Bbb=msbm10

\def\Z{\hbox{\Bbb Z}}

\numberwithin{equation}{section}

\DeclareMathOperator{\Aut}{Aut}

\newcommand\G{\Gamma }
\newcommand\D{{ \Delta }}
\newcommand\z{{ \zeta }}

        \def\<{{\langle}}
        \def\>{{\rangle}}

        \font\Aaa=msam10

\font\Aaa=msam10

\font\Bbb=msbm10

\def\Z{\hbox{\Bbb Z}}

\def\d{\delta}

\def\div{ \kern-.5pt\hbox{\big |} }
\def\ndiv{ {\not\kern-.5pt\hbox{\big |}\,} }
\def\ndivv{ {\not\kern+1.5pt\hbox{$\mid$}\,} }

\def\Aut{{\rm Aut}}

\def\col{\colon\!}

\def\B{^2\kern-.8pt B}
\def\G{^2\kern-.8pt G}
\def\EH{^2\kern-.8pt\hat  E}
\def\E{^2\kern-.8pt E}
\def\D{^3\kern-1pt D}
\def\FF{^2\kern-.8pt F}

\newdimen\refcodesize
\newbox\seriesbox
\def\proof{\noindent {\bf Proof.~}}

\def\ZZ{{\Z}}

\DeclareRobustCommand{\SkipTocEntry}[4]{}
\begin{document}

\title[Bent functions generalizing Dillon's partial spread functions]
{Bent functions generalizing Dillon's partial spread functions}
 \thanks{This research was supported in part by a grant from the NSA}
 
    \author{William M. Kantor
     }
      \address{College of Computer and Information Science, Northeastern U., Boston, MA 02115
}
                   \email{kantor@uoregon.edu}

%\dedicatory{\bigskip   \bigskip}

\begin{abstract} 
 \vspace{-6pt}
 This note presents  generalizations  of  the partial spread bent functions introduced by Dillon, as well as the corresponding relative difference sets in nonabelian groups.
\end{abstract} 
%\date{\today}
\keywords{bent function, partial spread, relative difference set}

\maketitle
 
 %\vspace{-20pt}
\section{Introduction}
In his thesis, Dillon introduced bent functions obtained using partial spreads
of $\ZZ_2$-vector spaces -- and, more generally, sets of subgroups of a group  \cite{Di,Di2}.  Generalized versions of this notion have been obtained using   desarguesian spreads
(\cite{Ho}, \cite[Theorem~40]{CD}, \cite{KGKN},\,\cite[Theorem~2.5]{Ny}) and partial spreads \cite{LL}.
 The present note provides a generalization
 (\thmref{nonabelian})
 of Dillon's bent functions involving groups,
   using what amount to partial 
  spreads  and  proved using   
  elementary bookkeeping but no exponential sums.
 All of the preceding results are special cases.    
 We conclude with another
  special case (\thmref{main}) using a proof involving 
   exponential sums.
   
   Since bent functions produce    relative difference sets
   \cite[p.~6]{Sch},
   we  obtain large numbers of relative difference sets  in  groups that need not be 
   abelian (\corref{corollary}).

  \section{Partial spreads of groups}
\label{Definitions}
  Let $G$ and  $H$ be finite groups.
\begin{definition}
\label{Def groups}
A function
$f\col G\to H$ is \emph{bent}  if, whenever 
$1\ne  z  \in G$, $x\mapsto f(xz )f(x)^{-1}$ takes each value in $H$ equally often (i.\,e., $|G|/|H|$ times).
\end{definition} 

This definition is trivially equivalent to the requirement that 
$\{  (x,f(x))\mid x \in  G\}$ is  a
 difference set in   $G\times H$  
 relative to $1 \times H\,$ \cite[p.~6]{Sch}.  
Such a function $f$  is called
``perfect nonlinear'' in \cite{Po}.
 
 The following is our main result:
      
     % \newpage
\begin{theorem}
\label{nonabelian}
 In
a  group $G$  of order $(qN)^2,$ 
let $\Sigma$ be a set  of   $(q-1)N$ subgroups of order $qN$ any two of which intersect only in $1$.  Let $H$ be a group of order $q$.
Partition $\Sigma$ into $q-1$  subsets  $\Sigma_i$ of size $N $
$(i\in H\backslash \{1\});$  let $D_i :=\cup\Sigma_i\backslash\{1\}$
and  $D_1 : =G\backslash\cup_{i\ne 1} D_i$. Then the function 
$f\col G\to  H, $ defined by
$f(D_i) =i $  for all $i\in  H,$
is bent. 
 
 \end{theorem}

%\newpage
\proof
It is crucial here  that (*)
  {\em $\{D_i  \mid  i\in H  \}$ partitions $G$.}
 We need to show that there are exactly $qN^2$  solutions $x$ to the equation $f(xz)f(x)^{-1}=b$ whenever $1\ne z\in G$
 and  $b\in  H$.  
  Thus, for each $z$ and $b$ we need to determine
 \begin{equation}
 \label{goal}
 \sum _{c,d}%
   |(D_cz^{-1})\cap D_d|  \mbox{ \ where 
 $  c,d\in   H$\,  satisfy\,
 $cd^{-1}=b$}.
  \end{equation}
  
  \vspace{-2pt}
  Let $k$ be the unique element of $H$ such that
   $z\in D_k$.  We always \emph{assume that  $i$ and  $j$ denote elements of  $H\backslash\{1\}$}, but $k$ may be 1.
If $k\ne1$   let  $z\in \tilde X_k\in \Sigma_k$.   (The tilde is included     in order to distinguish among subgroups in the same set $\Sigma_i$
 when $k=i$:  if $k\ne i$   then we delete the tilde.)

We proceed in several steps.
\begin{equation}
\label{Di}
 |D_i|=N(qN-1),  \ 
  |D_1|=qN^2 + qN- N.
\end{equation}
The first of these follows from $|\Sigma_i|=N$, and the second  from (*).
\begin{equation}
\label{Dij}
|(D_iz^{-1})\cap D_j|=(N-\d_{ik} ) (N-\d_{jk}) \
\mbox{ if\, $i\ne j$}.
\end{equation}
For, if 
$X_i\in \Sigma_i, X_j\in \Sigma_j, $ then
the equation $z=x_j^{-1}x_i$ ($x_i\in X_i, \, x_j\in X_j $) 
has a (unique) solution
for  distinct $i,j,k$, 
and  no solution if   $ X_i=  \tilde X_k$  or   $ X_j=  \tilde X_k$.
\begin{equation}
\label{Dii}
|(D_iz^{-1})\cap D_i|=(N-\d_{ik} ) (N-\d_{ik}-1)
+(qN-2)\d_{ik}.
\end{equation}
For, 
if    $X_i,  X_i'\in \Sigma_i$, then
the equation $z= x_i'{}^{-1}x_i$  (${x_i\in X_i,  \, x_i'\in X_i' }$) 
has a unique solution
  precisely when 
 $z\notin 
X_i,  X_i'$ and $X_i\ne X_i'$;   $\,qN-2$ solutions 
when  $ X_i= X_i' =\tilde X_k$
(since we must have $x_i,x_i'\ne1$); and no solution otherwise.
\begin{equation}
\label{Di0}
|(D_iz^{-1})\cap D_1|=
 |(D_1 z ^{-1})\cap D_i|=
 (N+1-  \d_{1k})
  (N-  \d_{ik})  +\d_{ik}  
  .
\end{equation} 
For, 
$|(D_iz^{-1})\cap D_1|=  |D_i z^{-1}|-\sum_{j\ne i} |(D_iz^{-1})\cap D_j|
-|(D_iz^{-1})\cap D_i |$ by  (*).  Now use \eqref{Di}-\eqref{Dii}
and an elementary calculation.   
\begin{equation}
\label{D00}
|(D_1 z^{-1})\cap D_1|=
 (N +1-  \d_{1k}) (N-  \d_{1k})
 + \d_{1k} qN  . 
\end{equation} 
For, 
$|(D_1z^{-1})\cap D_1|=  |D_1|-\sum_{i} |(D_iz^{-1})\cap D_1|$ by  (*).
 Now use \eqref{Di} and \eqref{Di0}.    

Two equally elementary calculations using \eqref{Di}-\eqref{D00}   and $\sum_c\d_{ck}=1$ show that  \eqref{goal}  
equals $qN^2$ for all $b\in  H$  (considering
the cases  $b\ne1$ and $b=1$ separately).\qed
 
% \newpage
 \medskip

{\noindent \bf Remark  1.  Group structure.}  
Unfortunately, 
if $G$ is  not elementary abelian,
 when $q=2$ there are very few examples of groups  having sets
 $\Sigma$ meeting our requirements   \cite{Fro}, and   
 when $q>2$ there are
no examples   
  \cite[Theorems~3.3, 3.4]{Ju}.  Therefore,   the preceding theorem only deals with elementary abelian groups $G$ when $q>2$.

  On the other hand, it 
seems unexpected that 
\emph{the definition of $f$ does not require any special properties of the  group $H$}, so that    nonisomorphic groups $H$ 
of order $q$ 
 produce  bent functions  
$f\col G\to  H $
 using  the same  partition of
  $\Sigma$. 
   
\medskip

{\noindent \bf Remark 2.  How many relative difference sets?}  It is not uncommon to provide    constructions of combinatorial objects and assume that   many inequivalent objects arise if there are many choices made in the construction
 \cite{Di,Di2,CD,Ho,KGKN,LL,Ny}.  Proving that there are, indeed, many inequivalent objects is another matter, one that can be difficult.

We already noted that bent functions correspond to relative difference sets.  Relative difference sets
in     nonisomorphic groups  $G\times H$  are clearly inequivalent.  There are certainly many examples showing that the structure of $H$ is usually involved
in such a relative difference set, unlike in the theorem.  

If $|H|= p^s$ with $p$ prime, then the number of nonisomorphic groups $H$ is large:
$p^ {(2/27)s^3+O(s^{8/3})}$ \cite{BNV}.
The number of inequivalent possibilities for the set $\Sigma$ is far larger \cite{Ka}.
However, cruder  estimates than in \cite{Ka}  already give   information 
in groups that need not be abeian (see   \cite[p.~6]{Sch}  for the parameters of a relative difference set): 
\begin{corollary}
\label{corollary}
For integers $m \ge s\ge1$ and a prime $p,$
let $G$ be an elementary abelian $p$-group of order $p^{2m}$
and let $H$ be any group of order $p^s$.
Then there are more than $p^{p^{m-1}-9m^2}$ pairwise inequivalent
 $(p^{2m}, p^{s}, p^{2m} ,p^{2m-s})$-difference sets in $G\times H$ relative to 
$1\times H$.
\end{corollary}

\medskip
\proof
``Inequivalence'' means ``in different $\Aut(G\times H)$-orbits''.
Clearly  $|\Aut(G\times H)|<p^{(2m+s)^2} \le p^{(3m)^2}$. The number of possible sets 
$\Sigma$ inside a desarguesian spread of $G$ is   
$% 
 \left(\begin{matrix}
p^m+1
\\
p^{m-1}
\end{matrix}\right) \ge p^{p^{m-1}}\!,$
producing more than
$p^{p^{m-1} -9m^2}$   inequivalent relative difference sets.  \qed

\medskip Note that the above estimate did not even
take  into account the many ways to partition a given choice $\Sigma$.
Far more bent functions and relative difference sets are obtained
from the Maiorana-McFarland bent functions 
(cf. \cite{Di2}, \cite[p.~51]{KSW} and 
\cite[Theorem~39]{CD})
using similar simple estimates, but those only use  elementary abelian
groups $G\times H$.

\medskip
{\noindent \bf Remark  3.  Association schemes. }
In the   notation of the theorem,  the sets 
$D_i$ produce   an association scheme,
obtained by partitioning $G\times G$ into the sets 
$\{(x,x ) \mid x\in G\}$ and 
$\{(x,y)\in G\times G \mid 1\ne xy^{-1}\in D_i\}$, $i\in H$
(compare \cite[esp.~p.~114]{vDM}). 
%\vspace{-10pt}

 \section{Vector spaces}
 We now turn to  a different type of proof of  a special case of \thmref{nonabelian},
 and another brief  discussion of   the number of different bent functions obtained.
 
  Let $V$ be a finite vector space over a   field $K$ of characteristic $p$.
There are two equivalent definitions of bent functions 
$V\to K$    \cite[Theorem~2]{Amb},  \cite[Theorem~2.3]{Ny}.
The first is  a special case of Definition~\ref{Def groups}:
\begin{definition} (Combinatorial definition.)
\label{Def 1}
$f\col V\to K$ is \emph{bent}  if, whenever 
$0\ne  z  \in V$, $v\mapsto f(v+ z )-f(v)$ takes each value in $K$ equally often (i.e., $|V|/|K|$ times).
\end{definition}

Let $\z    $ denote a primitive complex $p$th root of 1.  Fix a 
   nonzero linear functional
  $T\col K\to \Z _p$,  as well as a basis and hence a dot product for the $K$-space $V$.  For $f\col V\to  K$  and $  k   \in K$,
  write   $f_   k   (v):=\z    ^{T(   k   f(v) )} $  and
$$ 
\hat f( u ):=\sum_{v\in V}\z    ^{T( u\cdot v+f(v ))}
, ~ u\in V.
\vspace{-4pt}
$$ 
 \begin{definition}  (Fourier   definition.)
  \label  {Fourier   definition.}
$f\col V\to  K$ is \emph{bent}  if $|\hat f_ k  (u)|=|V|^{1/2}$
for all $ k \in K^*,  u \in V$. (This notion is independent of the choice of dot product and  $T$.)
\end{definition}
A function is   {\em balanced} if each member of the codomain occurs as a value equally often (compare 
Definitions~\ref{Def groups} and \ref{Def 1}.  This amounts to  a labelled partition of~the domain into sets of equal size, where the number of parts is the size of the codomain. 

 A finite    {\em prequasifield}    $(F,+,*)$ of characteristic $p$ consists of a finite vector space $F$ over $\Z_p$, together with a binary  operation $*$ on $F$ such that $a*(x+y)=a*x+a*y$ and $z\mapsto a*z-b*z$ is bijective for all $x,y,a,b\in F$, $a\ne b$.  The associated 
 {\em spread} 
  consists of the following $|F|+1$ subspaces of $F\oplus F$:  $x=0,$ and all 
 $y=m*x$ for $m\in F$;
 note the similarity to the situation in \thmref{nonabelian}.
The  associated {\em kernel}   is the field consisting of all additive maps 
$k\col F\to F$ such that $m*(kx)=k(m*x)$ for all $m,x\in F$.
Both the prequasifield and the spread determine a finite affine plane  \cite[p.~220]{De} that will not be needed here;
nor will the fact that
the same spread  can arise from   many non-isomorphic prequasifields.
 
We use  \defref {Fourier   definition.}
in order to provide an entirely different type of proof of  the following special case of \thmref{nonabelian}:

\begin{theorem}
\label{main}
Let $(F,+,*)$ be a prequasifield of characteristic $p$ whose kernel contains the field $K$.   Fix a  $K$-basis of $F$ and hence  of 
$V:=F\oplus F,$
and equip $ F$ and $V$ with the corresponding  dot products.   Let
$g\col  F\to K$ be any balanced function. Then $f\col V\to K$ is a bent function$,$  where
$f(0,y):=g(0) ,$
and
 $f(x,y):=g(m) $  with $y=m*x$ for 
 a unique $m\in F$ when $x\ne0$.
 \end{theorem}

\proof (Compare \cite[Theorem~40]{CD}.)~Clearly, 
$ 
%\displaystyle
 \hat f_   k   ( a , b )=
\sum_{x,y}\z    ^{T(( a , b )\cdot (x,y) +  k     f(x,y) )}.
$ 
 If $m\in F$  let $L_m\col F\to F$ be the $K$-linear map defined by $L_m(x)=m*x$; and let $L_m^t$ be its transpose, so  that
$b\cdot L_m(x)=L_m^t(b)\cdot x$ for all $b,x\in F$.
If $ x\in F^*$ and $y\in F$  then below we will write $(x,y)=(x,m*x)$
for a unique $m\in F$. 
For each $ k  \in  K^*$ and  $(a,b ) \in V$, 
 \begin{eqnarray*}
\hat f_   k   ( a , b ) \hspace{-4pt}
&=\hspace{-4pt}&\hspace{-4pt}\sum_{\scriptstyle 
x\in  F^*, % 
m\in F}\z    ^{T(( a , b )\cdot (x,m*x)+    k   f( x,m*x) )}
+\sum_{y\in F} \z    ^{  T(b \cdot y+  k   f(0,y)) }
\\
&=\hspace{-4pt}&\hspace{-4pt}\sum_{\scriptstyle 
x\in  F^*, % 
m\in F}\z    ^{T( a \cdot x+ b \cdot (m*x)) }\z    ^{T(   k   g(m ) )}
+\sum_{y\in F} \z    ^{  T(b \cdot y+    k   g(0) ) }
\\
&=\hspace{-4pt}&
\sum_{ m\in F }
\z    ^{ T(   k   g(m ) )}
\sum_{  x\in F^*}\z    ^{T( a \cdot x+ b \cdot L_m(x) )}
+\sum_{y\in F} \z    ^{  T( b \cdot y+   k   g(0) )  }
\\
&=\hspace{-4pt}&
\sum_{ m\in F }
\z    ^{ T(  k   g(m ))}
\sum_{ x \in F}\z    ^{T([ a  +L_m^t( b )] \cdot x) } 
-  \sum_{m\in F} \z    ^{ T(    k   g(m )) }
+\sum_{y\in F} \z    ^{ T(b \cdot y+  k   g(0) )}
\\
&=\hspace{-4pt}&
\sum_{ m\in F }
\z    ^{T(   k   g(m ) )}
\sum_{ x \in F}\z    ^{T([ a  +L_m^t( b )] \cdot x) } 
+\sum_{y\in F} \z    ^{ T(  b \cdot y+   k   g(0) ) };
\end{eqnarray*}
here    
$\sum_{ m }
\z    ^{  T( k   g(m ) )}=0$ since  
$\sum_0^{p-1}\z^j=0$ and
$m\mapsto   T( k   g(m ))$~is~balanced.

The transformations $L_m,m\in F$, have the property that the difference of any two is nonsingular; hence the same is
true of their transposes $L_m^t,m\in F$,
so that $m\mapsto L_m^t(b)$ is 1-1 and hence onto
if $ b\ne0$.  
Given $ b\ne0$ and $ a $ it follows that there is a unique 
$\tilde m\in F$ such that 
$ a  +L_{\tilde m}^t( b )=0$.  For that $\tilde m$ and  each $m'\ne \tilde m$, 
we have 
$\sum_{ x \in F}\z    ^{T([ a  +L_{m'}^t( b )] \cdot x )} =0$
since
$x\mapsto T( [ a  +L_{m'}^t( b )] \cdot x )$ is balanced.
Since $ y\mapsto T(  b \cdot y+   k   g(0) ) $ is also balanced, 
 $\hat f_   k   ( a , b )=  \z    ^{T(kg( \tilde m))}\sum_{x\in F}\z    ^0+0
= \z    ^{T(kg(\tilde m ) )}|F| $ has absolute value $|F|$.
 
 Finally, when  $b=0$     we find that 
$\hat f_   k   (a,0 )= \sum_{ m\in F }
\z    ^{T(   k   g(m ) )}
\sum_{ x \in F}\z    ^{T( a  \cdot x) } 
+\sum_{y\in F} \z    ^{ T(  k   g(0) ) }
=
%\break
0\sum_{ x \in F}\z    ^{T( a  \cdot x) }+ \z    ^{T(   k   g(0 ) )} |F| $ has absolute value $|F|$.\qed

\smallskip\smallskip\smallskip
{\noindent \bf Remark 4.} 
By \cite{Ka},  using  subsets $\Sigma$ of a desarguesian spread 
in \thmref{main} 
(so the quasifield is just a field)   
produces at least 
$ \left(\begin{matrix}
q^m+1
\\
q^{m-1}
\end{matrix}\right) \!
\Big/2(q^m+1)q^m(q^m-1)^2\log_p q^m$
pairwise affinely-inequivalent bent functions on $V$
(compare Remark~2).
Although there are many many different types of nondesarguesian spreads known, there are not enough known to change the preceding estimate 
significantly.

%\vspace{-8pt}

\end{document}